\documentclass{amsart}
 \usepackage{epsfig}
\usepackage{graphicx}
\usepackage{amscd}
\usepackage{amsmath}
\usepackage{amsxtra}
\usepackage{amsfonts}
\usepackage{amssymb}

\newtheorem{theorem}{Theorem}[section]
\newtheorem{corollary}[theorem]{Corollary}
\newtheorem{lemma}[theorem]{Lemma}
\newtheorem{proposition}[theorem]{Proposition}

\newtheorem{conjecture}[theorem]{Conjecture}
\theoremstyle{definition}
\newtheorem{definition}[theorem]{Definition}
\newtheorem{remark}[theorem]{Remark}

\newtheorem*{question}{Question}

\theoremstyle{remark}

\renewcommand{\theclaim}{\textup{\theclaim}}

\newtheorem*{acknowledgements}{Acknowledgements}

\numberwithin{equation}{section}

\def\openone

{\mathchoice

{\hbox{\upshape \small1\kern-3.3pt\normalsize1}}

{\hbox{\upshape \small1\kern-3.3pt\normalsize1}}

{\hbox{\upshape \tiny1\kern-2.3pt\SMALL1}}

{\hbox{\upshape \Tiny1\kern-2pt\tiny1}}}

\makeatletter

\newbox\ipbox

\newcommand{\ip}[2]{\left\langle #1\, , \,#2\right\rangle}
\newcommand{\diracb}[1]{\left\langle #1\mathrel{\mathchoice

{\setbox\ipbox=\hbox{$\displaystyle \left\langle\mathstrut
#1\right.$}

\vrule height\ht\ipbox width0.25pt depth\dp\ipbox}

{\setbox\ipbox=\hbox{$\textstyle \left\langle\mathstrut
#1\right.$}

\vrule height\ht\ipbox width0.25pt depth\dp\ipbox}

{\setbox\ipbox=\hbox{$\scriptstyle \left\langle\mathstrut
#1\right.$}

\vrule height\ht\ipbox width0.25pt depth\dp\ipbox}

{\setbox\ipbox=\hbox{$\scriptscriptstyle \left\langle\mathstrut
#1\right.$}

\vrule height\ht\ipbox width0.25pt depth\dp\ipbox}

}\right. }

\newcommand{\dirack}[1]{\left. \mathrel{\mathchoice

{\setbox\ipbox=\hbox{$\displaystyle \left.\mathstrut
#1\right\rangle$}

\vrule height\ht\ipbox width0.25pt depth\dp\ipbox}

{\setbox\ipbox=\hbox{$\textstyle \left.\mathstrut
#1\right\rangle$}

\vrule height\ht\ipbox width0.25pt depth\dp\ipbox}

{\setbox\ipbox=\hbox{$\scriptstyle \left.\mathstrut
#1\right\rangle$}

\vrule height\ht\ipbox width0.25pt depth\dp\ipbox}

{\setbox\ipbox=\hbox{$\scriptscriptstyle \left.\mathstrut
#1\right\rangle$}

\vrule height\ht\ipbox width0.25pt depth\dp\ipbox}

} #1\right\rangle}

\newcommand{\cj}[1]{\overline{#1}}

\newcommand{\bz}{\mathbb{Z}}

\newcommand{\br}{\mathbb{R}}
\newcommand{\bc}{\mathbb{C}}
\newcommand{\bt}{\mathbb{T}}
\newcommand{\bn}{\mathbb{N}}

\def\blfootnote{\xdef\@thefnmark{}\@footnotetext}

\renewcommand{\mod}{\operatorname{mod}}

\input xy
\xyoption{all}
\usepackage{amssymb}

\newcommand{\supp}[1]{\operatorname*{supp} (#1)}

\def\F{\mathcal{F}}

\def\H{\mathcal{H}}

\def\-{^{-1}}

\def\U{\mathcal{U}}

\def\ty{\emptyset}

\begin{document}
\title[Orthogonal exponentials]{Orthogonal exponentials, translations, and Bohr
completions}
\author{Dorin Ervin Dutkay}
\blfootnote{}
\address{[Dorin Ervin Dutkay] University of Central Florida\\
	Department of Mathematics\\
	4000 Central Florida Blvd.\\
	P.O. Box 161364\\
	Orlando, FL 32816-1364\\
U.S.A.\\} \email{ddutkay@mail.ucf.edu}

\author{Deguang Han}
\address{[Deguang Han] University of Central Florida\\
	Department of Mathematics\\
	4000 Central Florida Blvd.\\
	P.O. Box 161364\\
	Orlando, FL 32816-1364\\
U.S.A.\\}
\email{dhan@pegasus.cc.ucf.edu}

\author{Palle E.T. Jorgensen}
\address{[Palle E.T. Jorgensen]University of Iowa\\
Department of Mathematics\\
14 MacLean Hall\\
Iowa City, IA 52242-1419\\}\email{jorgen@math.uiowa.edu}
\thanks{Work supported in part by a grant from the National Science Foundation.} 

\begin{abstract}

We are concerned with an harmonic analysis in Hilbert spaces
$L^2(\mu)$, where $\mu$  is  a probability measure on $\br^n$. The
unifying question is the presence of families of orthogonal
(complex) exponentials $e_\lambda(x) = \exp(2\pi i \lambda x)$ in
$L^2(\mu)$. This question in turn is connected to the existence of
a natural embedding of $L^2(\mu)$ into an $L^2$-space of Bohr
almost periodic functions on $\br^n$. In particular we explore
when $L^2(\mu)$ contains an orthogonal basis of $e_\lambda$
functions, for $\lambda$ in a suitable discrete subset in $\br^n$;
i.e, when the measure $\mu$ is spectral. We give a new characterization of finite spectral sets in terms of the existence of a group of local translation. We
also consider measures $\mu$ that arise as fixed points (in the sense of Hutchinson) of
iterated function systems (IFSs), and we specialize to the case
when the function system in the IFS consists of affine and
contractive mappings in $\br^n$. We show in this case that if
$\mu$ is then assumed spectral then its partitions induced by the
IFS at hand have zero overlap measured in $\mu$. This solves part
of the \L aba-Wang conjecture. As an application of the new
non-overlap result, we solve the spectral-pair problem for
Bernoulli convolutions advancing in this way a theorem of Ka-Sing
Lau. In addition we present a new perspective
on spectral measures and orthogonal Fourier exponentials via the
Bohr compactification.
\end{abstract}
\maketitle \tableofcontents
\section{Introduction}\label{intr}

We explore the following general Fourier duality for probability measures $\mu$ with support contained in Euclidean space $\br^n$.
For vectors $\lambda$ in $\br^n$, we set $e_\lambda:= \exp (2\pi i \lambda\cdot x)$, and we consider each $e_\lambda$ as an element in the Hilbert space $L^2(\mu)$.
If $\Lambda$ is a subset in $\br^n$, we set $E(\Lambda):= \{ e_\lambda\, |\, \lambda \in \Lambda\}$.
If $E(\Lambda)$ is an orthonormal basis (ONB) in $L^2(\mu)$ we say that  the two $(\Lambda, \mu)$ form  a {\it spectral pair}. Spectral pairs have received recent attention in for example \cite{Fu74, JP92, JP93, JP94, JP98a, JP98b}; and we draw on results and motivation from these papers.

 In this paper, we give a new characterization of finite spectral sets.
We solve part of the \L aba-Wang conjecture \cite{LW02}. To do it, we offer more
powerful theorem in the context of affine IFSs allowing non-overlap (details
below); and we prove some new relations which are satisfied by affine
spectral measures (Proposition \ref{pra1}). As an application,  we solve the
spectral-pair problem for Bernoulli convolutions advancing in this way a
theorem of Ka-Sing Lau (from Adv. Math. 2007 \cite{HL08});  and we give a new
perspective on spectral measures and orthogonal Fourier exponentials via the
Bohr compactification.

   We explore the following questions/problems for orthogonal exponentials
in $L^2(\mu)$:
 (1) What geometric properties of $\mu$ are implied by the presence of
orthogonal exponentials? (2) We give conditions on the Fourier transform
$\hat\mu$ of a probability measure $\mu$ on $\br^n$, and subsets $\Lambda$ in $\br^n$
which are equivalent to the subset $\Lambda$ forming an orthogonal family of
exponentials in $L^2(\mu)$; and (3) a condition on $\hat\mu$ which characterizes
the case of such maximal orthogonal families.
    (4) In the case when $L^2(\mu)$ has an orthogonal basis of exponentials
(ONB), we show that $\mu$ itself is determined by a family of local
translations (defined in the paper). This is accomplished with the use of a
family of unitary representations of the additive group $\br^n$, and Bohr's
theory of almost periodic functions.
  (5) In the special case when $\mu$ is an IFS-measure in the sense of
Hutchinson, we show that if $\mu$ has an ONB of exponentials, i.e., if $\mu$ is
spectral, then the subdivided parts of $\mu$ must have non-overlap. (6) We
then use this to prove that if $\mu$ is an infinite convolution-Bernoulli
measure with scale $\lambda$, then all the cases of measures $\mu_\lambda$, for
$\lambda > 1/2$ are measures of non-spectral type. (7) Finally we offer
detailed information about $\mu_\lambda$  in the case $\lambda = 3/4$.

   We will recall the Bohr-Besicovitch $L^2$-almost periodic compactification $G$ of $\br^n$; \cite{Bes32, BeBo31, Bes32, Bes35}. Bohr's group $G$ is an almost periodic completion, a compact group, so with normalized Haar measure; and its discrete dual group of characters is the group $\br^n$ with the discrete topology.

    This means that for each $\lambda$ in $\br^n$, $e_\lambda$ is then viewed as a character on $G$.

Let $\mu$ be a probability measure on $\br^n$, Borel. Suppose, there is a subset $\Lambda$ in $\br^n$ such that $E(\Lambda)$ is an ONB in $L^2(\mu)$( so $(\Lambda, \mu)$ is a spectral pair); then it follows that $L^2(\mu)$ embeds isometrically into $L^2(G)$, with the isometry determined by sending $e_\lambda$ in $L^2(\mu)$ into $\tilde e_\lambda$ in $L^2(G)$, see Theorem \ref{th5.3}.

The converse is true too \cite{JP98a, JP98b}.

   Let $\Omega$ be a subset in $\br^n$ of finite positive Lebesgue measure. Specializing now to $\mu:=$ Lebesgue in $\br^n$ restricted to some $\Omega$ in $\br^n$, Fuglede considered in 1974 \cite{Fu74} these measures, and he suggested (conjectured?) that a given $\Omega$ has the spectral property if and only if it tiles $\br^n$ under translations by points in $\br^n$. Making the connection to Bohr's theory of almost periodic functions, one now sees an intuition behind Fuglede's conjecture.

    If the embedding result sketched above were periodic, as opposed to almost periodic, then we would have the truth of the Fuglede conjecture. As it turns out it was negative (Tao \cite{Ta04}). And with hindsight we note that a negative answer was reasonably to to be expected. As for a positive result, see however \cite{IKT03}.

\begin{definition}\label{defi1}
Let $e_\lambda(x):=e^{2\pi i\lambda\cdot x}$, $x,\lambda\in\br^n$. Let $\mu$ be a Borel probability measure on $\br^n$. We say that $\mu$ is a {\it spectral measure} if there exists a subset $\Lambda$ of $\br^n$ such that the restrictions of the functions $e_\lambda$, $\lambda\in\Lambda$ form an orthogonal basis for $L^2(\mu)$. In this case $\Lambda$ is called a {\it spectrum} for the measure $\mu$.
\end{definition}

  The paper is organized as follows: In section \ref{seclt} we introduce the notion of ``local
translations'' associated to a spectral measure. In section \ref{secat} we characterize all the
atomic spectral measures in terms of the existence of a group of local translations, and we return to atoms again in Proposition \ref{prnoa}.
It turns out that if some spectral measure $\mu$ with spectrum $\Lambda$ has just
one atom, then $\Lambda$ is necessarily finite, and the situation is covered
by our theorem in Section \ref{sec3}.  This is then applied in section \ref{sec3} where we show that if a
spectral pair $(\mu, \Lambda)$ arises as a Hutchinson measure associated with
an affine IFS, then (Theorem \ref{thnoo}) $\mu$ is a ``no-overlap'' IFS measure.
Finally in section \ref{sec4} we construct embeddings of spectral measures into the Bohr compactification and show that they intertwine the local translations. 

\section{The group of local translations}\label{seclt}

   Historically spectral pairs arose \cite{Fu74} in the study of domains $\Omega$
in $\br^n$ and consideration of the partial derivative operators in $L^2(\Omega)$
for a suitable domain $\Omega$ in $\br^n$; so the vector fields in the $n$
coordinate directions, defined on compactly supported $C^1$ functions in
$\Omega$. This problem lends itself naturally to the consideration of local
translations (see also \cite{JP93, JP99, Jor82}). In this section we extend this
idea of local translations to arbitrary spectral pairs, and we derive some
consequences.

\begin{definition}\label{defft}\cite{JP99}. 
Let $\mu$ be a spectral probability measure on $\br^n$, with spectrum $\Lambda$ a subset of $\br^n$.

Define the Fourier transform $\F:L^2(\mu)\rightarrow l^2(\Lambda)$ by 
$$(\F f)(\lambda)=\ip{f}{e_\lambda},\quad(f\in L^2(\mu),\lambda\in\Lambda).$$
Then $\F$ is unitary and  
$$\F^{-1}(c_\lambda)_\lambda=\sum_{\lambda\in\Lambda}c_\lambda e_\lambda.$$

Define the group of transformations $(U(t))_{t\in\br^n}$ on $L^2(\mu)$ by
$$U(t)f=\F^{-1}((e_t(\lambda)\F f(\lambda))_\lambda)=\sum_{\lambda} [e^{2\pi i t\cdot\lambda}\ip{f}{e_\lambda}] e_\lambda.$$
The convergence of the sum is in $L^2(\mu)$.
This means that in the ``Fourier domain'', $\hat U(t):=\F U(t) \F^{-1}$ is just multiplication by the sequence $(e^{2\pi i t\cdot\lambda})_\lambda$. We call $(U(t))_{t\in\br^n}$ the {\it group of local translations}.

Note also that 
\begin{equation}
	U(t)e_\lambda=e_\lambda(t)e_\lambda,\quad(t\in \br^n,\lambda\in\Lambda).
	\label{eqeig}
\end{equation}
Note that $U(t)$ depends on the spectrum $\Lambda$. 
\end{definition}

The reason for the name ``local translations'' is given in the following Proposition.
\begin{proposition}\label{prjp}\cite{JP99}
Let $\mu$ be a spectral measure on $\br^n$ and with spectrum $\Lambda$, and let $(U(t))_{t\in\br^n}$ be its group of local translations. Assume that $\mu$ is compactly supported.
Suppose $O\subset\br^n$ is measurable, $t\in\br^n$, and $O, O+t\subset\supp\mu$, where by $\supp\mu$ we denote the {\it support }of the measure $\mu$, i.e., the smallest compact set $X$ with $\mu(X)=1$. Then 
$$(U(t)f)(x)=f(x+t)$$
for a.e. $x\in O$ and every $f\in L^2(\mu)$. Moreover
$$\mu(O+t)=\mu(O).$$
\end{proposition}

\begin{corollary}\label{cor2.3}
If $O,O+t\subset\supp\mu$ then 
$$U(t)\chi_{O+t}=\chi_O.$$
\end{corollary}

\begin{proof}
From Proposition \ref{prjp}, we have that $U(t)\chi_{O+t}(x)=\chi_{O+t}(x+t)=1=\chi_O(x)$ for $\mu$-a.e. $x\in O$. But
$U(t)$ is unitary so 
$$\mu(O+t)=\mu(O)=\int_O1\,d\mu\leq\int_O1\,d\mu+\int_{\br\setminus O}|U(t)\chi_{O+t}|^2\,d\mu=$$$$\int |U(t)\chi_{O+t}|^2\,d\mu=\|U(t)\chi_{O+t}\|^2=\|\chi_{O+t}\|^2=\mu(O+t).$$
Then $U(t)\chi_{O+t}(x)=0$ for $\mu$-a.e. $x\in\br^n\setminus O$.
This implies the corollary.
\end{proof}

We propose the following question on a possible characterization of spectral measures. As we will see in section \ref{secat}, the question has a positive answer for atomic measures:

\begin{question}
Let $\mu$ be a probability measure on $\br^n$. Suppose there is a (strongly continuous) group of unitary transformations $(U(t))_{t\in\br^n}$ on $L^2(\mu)$ such that for every measurable set $O$, and $t\in\br^n$ with $O,O+t\subset\supp\mu$
$$U(t)\chi_{O+t}=\chi_O.$$
Is $\mu$ is a spectral measure?
\end{question}

\subsection{Atomic spectral measures}\label{secat}

   In this section we find the spectral pairs $(\mu, \Lambda)$  for which $\mu$
is a sum of Dirac masses, and we give a characterization of such spectral measures in terms of the existence of a group of local translations.

\begin{definition}
Let $A$ be a finite subset of $\br^n$. $N:=\#A$. We say that $A$ is a {\it spectral} set if the normalized counting measure $\delta_A$ on $A$ is a spectral measure, $\delta_A:=\frac{1}{N}\sum_{a\in A}\delta_a$. A spectrum for $A$ is a spectrum for the measure $\delta_A$. We denote by $L^2(A):=L^2(\delta_A)$.
\end{definition}

\begin{theorem}\label{thi1}
Let $A$ be a finite subset of $\br^n$. The following affirmations are equivalent:
\begin{enumerate}
\item
The set $A$ is spectral.
\item
There exists a continuous group of unitary operators $(U(t))_{t\in\br^n}$ on $L^2(A)$, i.e., $U(t+s)=U(t)U(s)$, $t,s\in\br^n$ such that 
\begin{equation}
	U(a-a')\chi_a=\chi_{a'}\quad(a,a'\in A),
	\label{eqthi1_1}
\end{equation}
where 
$$\chi_a(x)=\left\{\begin{array}{cc}
1,&x=a\\
0,&x\in A\setminus\{a\}.\end{array}\right.$$
\end{enumerate}
\end{theorem}

\begin{proof}
(i)$\Rightarrow$(ii). Follows from Corollary \ref{cor2.3}.

(ii)$\Rightarrow$(i). By Stone's theorem (the multivariable version, Stone-Naimark-Ambrose-Godement \cite{Am44}), there exist commuting self-adjoint operators $G_j$ on $L^2(A)$ such that $U(0,\dots,t_j,\dots,0)=e^{2\pi i t_jG_j}$ for all $t_j\in\br$ and all $j\in\{1,\dots,n\}$. Let $\{v_1,\dots,v_N\}$ be an orthonormal basis of eigenvectors common to $G_1,\dots,G_n$, $G_jv_i=\lambda_j^{(i)}v_i$ for $i\in\{1,\dots,N\}, j\in\{1,\dots, n\}$. Let $P_i$ be the orthogonal projection onto $v_i$. The projections are orthogonal. We have 
$$U(0,\dots,t_j,\dots,0)=\sum_{i=1}^Ne^{2\pi i\lambda_j^{(i)}t_j}P_i,\quad(t_j\in\br).$$
Therefore
$$U(t_1,\dots,t_n)=\sum_{i=1}^Ne^{2\pi i(\lambda_1^{(i)}t_1+\dots+\lambda_n^{(i)}t_n)}P_i$$

Let $\lambda^{(i)}:=(\lambda_1^{(i)},\dots,\lambda_n^{(i)})$, $\Lambda:=\{\lambda^{(1)},\dots,\lambda^{(N)}\}$. We relabel $P_{\lambda^{(i)}}:=P_i$ and $v_{\lambda^{(i)}}:=v_i$. 
We have

%
%
$$U(t)=\sum_{\lambda\in\Lambda} e^{2\pi i\lambda\cdot t}P_\lambda,\quad(t\in\br^n).$$
Using \eqref{eqthi1_1}, we have for $a,a'\in A$:
$$\sum_{\lambda}e^{2\pi i\lambda\cdot(a-a')}P_\lambda(\chi_a)=\chi_{a'}.$$
Applying $P_{\lambda'}$ to both sides we get $e^{2\pi i\lambda'\cdot a}P_{\lambda'}(\chi_a)=e^{2\pi i\lambda'\cdot a'}P_{\lambda'}(\chi_{a'})$. This implies that, for each $\lambda\in\Lambda$, there exists some constant $c(\lambda)\in\bc$ such that $P_\lambda(\chi_a)=e^{-2\pi i\lambda\cdot a}c(\lambda)v_\lambda$, for all $a\in A$. Then
$$\chi_a=\sum_{\lambda}P_\lambda(\chi_a)=\sum_\lambda e^{-2\pi i\lambda\cdot a}c(\lambda)v_\lambda.$$
Since $\chi_a$ and $\chi_{a'}$ are orthogonal if $a\neq a'$, we obtain
\begin{equation}
\sum_\lambda e^{2\pi i(a-a')\cdot\lambda} |c(\lambda)|^2=0.	
	\label{eqthi1_3}
\end{equation}

But then this shows that the measure $\mu_{c,\Lambda}:=\sum_{\lambda}|c(\lambda)|^2\delta_\lambda$ has the orthogonal exponentials $\{e_a\,|\,a\in A\}$, and since we also have $\#\Lambda\leq N=\#A$, it follows that $\mu_{c,\Lambda}$ is spectral with spectrum $A$, and $\#\Lambda=\#A=N$. Then, using \cite[Theorem 1.2]{LW06}, or Proposition \ref{prjp}, we obtain $|c(\lambda)|=$constant. Then \eqref{eqthi1_3} implies that the matrix $(e^{2\pi ia\lambda})_{a\in A,\lambda\in\Lambda}$ has orthogonal rows, so it has orthogonal columns, so $\Lambda$ is a spectrum for $A$.
\end{proof}

In the next theorem we give a different characterization of atomic spectral measures in terms of existence of functions which resemble the reproducing kernels for Hilbert spaces.

\begin{theorem}\label{thi2}
Let $A$ be a subset of $\br^n$. The following affirmations are equivalent:

\begin{enumerate}
\item
$A$ is a spectral set.
\item There exists a continuous function $c:\br^n\rightarrow\bc$ with the following properties:
\begin{equation}
	c(a'-a)=\delta_{aa'},\quad(a,a'\in A).
	\label{eqt2.1}
\end{equation}
\begin{equation}
	c(-t)=\cj{c(t)},\quad(t\in\br^n).
	\label{eqt2.2}
\end{equation}
\begin{equation}
	c(u_1-u_2)=\sum_{a\in A}c(u_1+a)\cj{c(u_2+a)},\quad(u_1,u_2\in\br^n).
	\label{eqt2.3}
\end{equation}
\end{enumerate}
\end{theorem}

\begin{proof}
(i)$\Rightarrow$(ii). We use Theorem \ref{thi1} and its proof. Let
\begin{equation}
	c(t)=\frac1N\sum_{\lambda\in\Lambda}e_\lambda(t),\quad(t\in \br).
	\label{eqt2.4}
\end{equation}

Then for $t\in\br^n$, $a,a'\in A$.
$$(U(t)(\chi_a))(a')=N\ip{U(t)\chi_a}{\chi_{a'}}=N\ip{U(t)\left(\sum_{\lambda\in\Lambda}\ip{\chi_a}{e_\lambda}e_\lambda\right)}{\chi_a'}=
N\sum_{\lambda\in\Lambda}\ip{\left(\frac1N e_{-\lambda}(a)e_t(\lambda)\right)e_\lambda}{\chi_{a'}}$$
Thus 
\begin{equation}
	(U(t)(\chi_a)(a')=\frac{1}{N}\sum_{\lambda\in\Lambda}e_\lambda(t+a'-a),\quad(t\in\br^n,a,a'\in A).
	\label{eqthi1_2}
\end{equation}
Therefore
\begin{equation}
	U(t)_{aa'}=N\ip{U(t)\chi_a}{\chi_{a'}}=c(t+a'-a),\quad(t\in\br^n,a,a'\in A).
	\label{eqt2.5}
\end{equation}
We have:
$c(a'-a)=U(0)_{aa'}=\delta_{aa'}$, since $U(0)=I_N$. Since $U(-t)=U(t)^*$, we have 
$$c(-t)=U(-t)_{aa}=\cj{U(t)_{aa}}=\cj{c(t)}.$$
Since $U(t_1+t_2)=U(t_1)U(t_2)$, we obtain 
$$c(t_1+t_2+a''-a)=U(t_1+t_2)_{a,a''}=\sum_{a'\in A}U(t_1)_{aa'}U(t_2)_{a'a''}=\sum_{a'}c(t_1+a'-a)c(t_2+a''-a').$$
Changing the variable to $u_1=t_1-a$ and $u_2=-t_2+a''$ we obtain \eqref{eqt2.3}.

(ii)$\Rightarrow$(i) 
Define $U(t)_{aa'}=c(t+a'-a)$ for all $t\in\br^n,a,a'\in A$.
Then for $a,a''\in A$, using \eqref{eqt2.2},
$$\sum_{a'}U(t_1)_{aa'}U(t_2)_{a'a''}=\sum_{a'}c(t_1+a'-a)c(t_2+a''-a')\stackrel{\mbox{\eqref{eqt2.3}}}{=}c(t_1+t_2+a''-a)=U(t_1+t_2)_{aa''},$$
so $U(t_1)U(t_2)=U(t_1+t_2)$.

Then $U(-t)_{aa'}=c(-t+a'-a)=\cj{c(t-a'+a)}=\cj{U(t)_{a'a}}$ so $U(-t)=U(t)^*$. 

Also $U(0)_{aa'}=c(a'-a)=\delta_{aa'}$ so $U(0)=I_N$. Finally $U(t)U(t)^*=U(t)U(-t)=U(0)=I_N$ so $U(t)$ is unitary.

\end{proof}

%
%

\section{Affine IFSs}\label{sec3}

   We exploit the idea of local translations for arbitrary spectral pairs,
and we reconsider a conjecture of \L aba and Wang. A main result is for the
particular spectral pairs $(\mu, \Lambda)$  where $\mu$ arises as a Hutchinson
measure associated with an affine and contractive IFS. In this case we show
(Theorem \ref{thnoo}) that $\mu$ must then be a ``no-overlap'' IFS measure.

Since iterated function systems(IFS) by definition involve a ``global''
set $X$ built by an iteration of copies of itself, i.e., there is an indexed and finite family of functions, say $(\tau_b)_{b\in B}$, and we are interested in possible overlap of the sets $\tau_b(X)$ as $b$ varies over the index set. There is an analogous selfsimilarity rule applied to specific measures $\mu$ supported on $X$. When the system $(\tau_b )_{b\in B}$ is given, the pair $\mu, X$ satisfies the so called Hutchinson condition \cite{Hut81}.

    More precisely, the question of ``overlap'' consists of deciding, for a given IFS, whether $\mu(\tau_b(X) \cap \tau_{b'}(X))$ is positive, when $b$ and $b'$ are distinct. So the possible overlap is measured in the same selfsimilar measure $\mu$, not in some different measure.

     Indeed, often there is clearly some overlap as the index point $b$ varies, but the question is the size of the overlap measured by  $\mu$. 
    In fact the study of IFSs naturally divides up into the two cases: (i) overlap may occur, and (ii) the case of ``non-overlap''; again referring to measurement of possible overlap in the specific selfsimilar measure $\mu$.

    As it turns out, for a particular IFS, it is often difficult to decide if we are in one or the other of the two cases; and the case of overlap, i.e., (i) is typically the more subtle of the two.

    Here we show that if an IFS is affine, and if it is further assumed to be spectral, then it follows that it is ``non-overlap''. This is of interest since it may be used in showing that a class of IFS measures $\mu$ (for example certain Bernoulli convolutions) are non-spectral measures.

     The earlier literature on the subject includes these papers: \cite{BJ99, DJ07a, DJ07b, DJ07c, HL08a, HL08b, JP92, JP94, JP95, JP96, JKS07a, JKS07b, LW05, LN07}. And a classical case where the overlap question arises is that of certain infinite Bernoulli convolutions; see details below, as well as the literature, for example \cite{FO03, FW03, Yu04}.

\begin{definition}\label{def3.1} Let $R$ be a $n\times n$ expansive real matrix, $S=R^t=$transpose of $R$, $B$ a finite subset of $\br^n$ containing 0, and let $N:=\#B$. Define $\tau_b$ by
$$\tau_b(x)=R^{-1}(x+b),\quad(x\in\br^n,b\in B).$$
By \cite{Hut81}, there is a unique compact set $X_B$ such that
$$X_B=\bigcup_{b\in B}\tau_b(X_B).$$
The set $X_B$ is called the {\it attractor} of the affine iterated function system (IFS) $(\tau_b)_{b\in B}$,

Given a list of probabilities $p_b\in(0,1)$, such that $\sum_{b\in B}p_b=1$, there is a unique probability measure $\mu_{B,p}$ such that for all continuous functions on $\br^n$,
\begin{equation}
	\int f\,\mu_{B.p}=\sum_{b\in B}p_b\int f\circ\tau_b\,d\mu_{B,p}.
	\label{eqh1}
\end{equation}
$\mu_{B,p}$ is called the {\it invariant measure} for $(\tau_b)_{b\in b}$ and $(p_b)_{b\in b}$. Moreover, the support of the measure $\mu_{B,p}$ is $\supp\mu=X_B$. When $p_b=1/N$ for all $b\in B$, we use the shorter notation $\mu_B$ or just $\mu$.
\end{definition}

\begin{conjecture}{\bf [\L aba-Wang \cite{LW02}]} \label{conlw}
Consider the dimension $n=1$. Let $\mu=\mu_{B,p}$ be the invariant measure associated with the IFS $(\tau_b)_{b\in B}$, with probability weights $(p_b)_{b\in B}$. Suppose that $\mu_{B,p}$ is a spectral measure. Then
\begin{enumerate}
\item[(a)] $R=\frac1M$ for some $M\in\bz$.
\item[(b)] $p_b=\frac1N$ for all $b\in B$.
\item[(c)] Suppose that $0\in B$. Then $B=\alpha D$ for some $\alpha\in\br$ and $D\subset\bz$. Furthermore, $D$ must be a complementing set $(\mod N)$, i.e., there exists a set $E\subset \bz$ such that $D\oplus E$ is a complete residue system $(\mod N)$. 
\end{enumerate}
\end{conjecture}

The following Proposition shows that part (b) of the \L aba-Wang conjecture is true when there is no overlap. 
\begin{proposition}
In the hypotheses of Conjecture \ref{conlw} assume there is no overlap, i.e. $\mu(\tau_b(X_B)\cap\tau_{b'}(X_B))=0$ for $b\neq b'$, where $X_B$ is the attractor of the IFS. Then $p_b=\frac1N$ for all $b\in B$. 

\end{proposition}

\begin{proof}

Fix $b_0\in B$. We have that $\tau_{b}(X)$ is contained in $\supp\mu=X_B$, for all $B$, and $$\tau_{b}(X_B)=R^{-1}(X_B+b)=R^{-1}(X_B+b_0)+R^{-1}(b-b_0)=\tau_{b_0}(X_B)+R^{-1}(b-b_0).$$
So $\tau_{b_0}(X_B)+R^{-1}(b-b_0)\subset\supp\mu=X_B$. Using Proposition \ref{prjp}, we get that $\mu(\tau_{b_0}(X_B))=\mu(\tau_{b}(X_B))$ for all $b\in B$. 
But, since there is no overlap, $\mu(\tau_b(X_B))=p_b$. Thus $p_b=$const$=\frac1N$.
\end{proof}

The next theorem shows that, in the case of equal probabilities, spectral measures have no overlap.

\begin{theorem}\label{thnoo}
Suppose the invariant measure $\mu=\mu_B$ (equal probabilities) is spectral. Then there is no overlap.
\end{theorem}

\begin{proof} We will need the following Lemma:
\begin{lemma}\label{prextm}
Suppose $\mu$ is a spectral measure with support $X\subset\br^n$. Let $\Gamma$ be a countable subgroup of $\br^n$. Then there exists a Borel measure $\tilde\mu$ on $\br^n$ such that
\begin{enumerate}
\item For all Borel subsets $A$ of $X+\Gamma$ and all $\gamma\in\Gamma$ we have $\tilde\mu(A)=\tilde\mu(A+\gamma)$.
\item The restriction of $\tilde\mu$ to $X$ is $\mu$. 
\end{enumerate}
\end{lemma}

\begin{proof}
Let $\Gamma:=\{\gamma_n\,|\,n\in\bn\}$ with $\gamma_1=0$. Let $A_n:=X+\gamma_n$. Let $B_1:=A_1=X$ and define $B_n$ as $B_n:=A_n\setminus\cup_{k=1}^{n-1}A_k$. The sets $B_n$ are disjoint, and $\cup_n B_n=\cup_n A_n=X+\Gamma$.

For a Borel set $A$ define
$$\tilde\mu(A)=\sum_{n=1}^\infty\mu[(A\cap B_n)-\gamma_n].$$
Clearly $\mu$ is a Borel measure, and the restriction of $\tilde\mu$ to $X=A_1$ is $\mu$.

Take now $A$ a Borel subset of $X+\gamma$, for some $\gamma\in\Gamma$. We have $\tilde\mu(A\cap B_n)=\mu[(A\cap B_n)-\gamma_n]$. Since $(A\cap B_n)-\gamma_n$ and $(A\cap B_n)-\gamma$ are contained in $X=\supp \mu$, it follows from the local translation invariance of $\mu$, that 
$\tilde\mu(A\cap B_n)=\mu[(A\cap B_n)-\gamma]$. This implies that $\tilde\mu(A)=\mu(A-\gamma)$.

It is enough to prove (i) for sets $A$ contained in one of the sets $B_n$. If $A\subset B_n\subset A_n$, and $\gamma\in\Gamma$, then 
$A+\gamma$ is contained in $A_n+\gamma=X+\gamma_n+\gamma=X+\gamma_m=A_m$ for some $m$. Therefore $$\tilde\mu(A+\gamma)=\mu(A+\gamma-\gamma_m)=\mu(A-\gamma_n)=\tilde\mu(A).$$
This proves (i).
\end{proof}

Returning to the proof of the theorem, let $X=\supp\mu$. We will prove first that for $b_1,b_2\in B$. 
\begin{equation}
	X\cap\bigcup_{b\neq b_1}(X+b-b_1)=X\cap\bigcup_{b\neq b_2}(X+b-b_2)\quad\mu\mbox{-a.e.}
	\label{eqov1}
\end{equation}
We reason by contradiction. Suppose there are $b_1\neq b_2\in B$ such that \eqref{eqov1} is not satisfied. Then there is a set $A\subset X$ such that
$\mu[A\cap(X+b_0-b_2)]>0$ for some $b_0\in B$, $b_0\neq b_2$, and $\mu[A\cap(X+b-b_1)]=0$ for all $b\neq b_1$ (or vice versa, in which case we switch between $b_1$ and $b_2$). 

Then using the local translation invariance of $\mu$ (Proposition \ref{prjp}), we have $\mu[(A+b_2-b_0)\cap X]>0$ and $\mu[(A+b_1-b)\cap X]=0$ for all $b\neq b_1$. 

From the invariance equation we have
$$\mu(\tau_{b_2}A)=\frac1N\sum_{b\in B}\mu(\tau_b^{-1}\tau_{b_2}A)=\frac{1}{N}\sum_{b\in B}\mu(A-b+b_2)\geq\frac1N\mu(A)+\frac{1}{N}\mu(A-b_0+b_2)>\frac{1}{N}\mu(A).$$
Also 
$$\mu(\tau_{b_1}A)=\frac1N\sum_{b\in B}\mu(A-b+b_1)=\frac1N\mu(A).$$
But, the translation invariance implies that $\mu(\tau_{b_1}A)=\mu(\tau_{b_2}A)$. This yields a contradiction.

Let $\Gamma$ be the subgroup of $\br^n$ generated by $B$. Consider the measure $\tilde\mu$ from Lemma \ref{prextm}. In what follows, all the inclusions will be $\tilde\mu$-a.e.

Fix $b_1\in B$. From \eqref{eqov1} we have that for all $b_2\neq b_1$, $\cup_{b\neq b_2}(X\cap (X+b-b_2))=\cup_{b\neq b_1}(X\cap(X+b-b_1))$, therefore $X\cap(X+b_1-b_2)\subset\cup_{b\neq b_1}(X\cap(X+b-b_1))$. 
Since $\tilde\mu$ is $\Gamma$-translation invariant we have 
$$(X+b_2-b_1)\cap X\subset \bigcup_{b\neq b_1}((X+b_2-b_1)\cap(X+b-b_1+b_2-b_1)),$$
so
$$\bigcup_{b_2\neq b_1}(X\cap (X+b_2-b_1))\subset \bigcup_{b_3\neq b_1,b_2\neq b_1}((X+b_2-b_1)\cap(X+b_3-b_1+b_2-b_1)).$$
By induction, using the translation invariance of $\tilde\mu$, we obtain that 
$$\bigcup_{b_2,\dots b_n\neq b_1}((X+\sum_{k=2}^{n-1}(b_k-b_1))\cap(X+\sum_{k=2}^n(b_k-b_1)))\subset
\bigcup_{b_2,\dots b_{n+1}\neq b_1}((X+\sum_{k=2}^{n}(b_k-b_1))\cap(X+\sum_{k=2}^{n+1}(b_k-b_1))).$$
Therefore, for all $n\geq2$,
\begin{equation}
X\cap\bigcup_{b_2\neq b_1}(X+b_2-b_1)\subset\bigcup_{b_2,\dots b_n\neq b_1}((X+\sum_{k=2}^{n-1}(b_k-b_1))\cap(X+\sum_{k=2}^n(b_k-b_1))).	
	\label{eqov2}
\end{equation}

The set on the left moves away from the initial set as $n$ increases, and this will give us a contradiction. We use the following lemma:
\begin{lemma}\label{lemleft}
Given a finite set of points $B$ in $\br^n$ there exists a unit vector $u\in\br^n$ and a $\delta>0$ such that there is a $b_1\in B$ with the property that $(b-b_1)\cdot u>\delta$ for all $b\in B$, $b\neq b_1$. 
\end{lemma}

\begin{proof}
Since the set $(B-B) \setminus \{0\}$ is finite, there exists a vector $u\in\br^n$ which is not orthogonal to any of the vectors $b_1-b_2$, $b_1\neq b_2$. So $b_1\cdot u\neq b_2\cdot u$ for all $b_1\neq b_2$. Take $b_1$ such that $b_1\cdot u=\min_{b\in B}b\cdot u$. Then $(b-b_1)\cdot u> 0$ for all $b\neq b_1$, and since $B$ is finite the lemma follows. 
\end{proof}
Take $u$ and $b_1\in B$ as in Lemma \ref{lemleft}. Then $(\sum_{k=2}^n(b_k-b_1))\cdot u>(n-1)\delta$ for all $n$, and all $b_2,\dots,b_n\neq b_1$. The map $X\times X\ni(x,y)\mapsto (x-y)\cdot u$ is bounded since $X$ is compact. Take $n$ such that $(n-1)\delta> \sup_{x,y\in X}(x-y)\cdot u$. Then $X\cap(X+\sum_{k=2}^{n-1}(b_k-b_1))=\ty$.
But then, using \eqref{eqov2}, since the left-hand side is contained in $X$, it follows that the inclusion can be realized only if the set on the left has measure $\tilde\mu$ zero. Then, using \eqref{eqov1} and since $\tilde \mu=\mu$ on $X$, it follows that $\mu(X\cap(X+b-b'))=0$ for all $b\neq b'$.
Then, since $\mu$ is supported on $X$, $\mu(X+b-b')=0$. 

Finally, using the invariance equation, we have 
$$\mu(\tau_b(X))=\frac{1}{N}\sum_{b'\in B}\mu(\tau_{b'}^{-1}(\tau_b(X))=\frac1N\sum_{b'\in B}\mu(X+b'-b)=\frac1N\mu(X)=\frac1N.$$
Since $\mu(\cup_b\tau_b(X))=\mu(X)=1$, it follows that there can be no overlap.

\end{proof}

The next proposition gives some positive evidence that part (c) of the \L aba-Wang conjecture might be true (after some modifications). As one can see in \eqref{eqa1_2}, we have a $\sum e^{2\pi i S^{-1}\lambda\cdot (b-b')}\dots=0$. If we could ignore the other terms (which we might be able to do since \eqref{eqa1_1} holds), then we would get some relation which is close to $B$ being a spectral set, which is a weaker form of (c) in the \L aba-Wang conjecture.
\begin{remark}\label{rem4.7}
We have to take into consideration the counterexample given in \cite{DJ08} with $R=4$, $B=\{0,1,8,9\}$ where the attractor is $[0,1]\cup[2,3]$ and $\mu$ is the Lebesgue measure on this set. This is a spectral measure, but $B$ is not complementing mod $4$. However $B$ is a spectral set.
\end{remark}
Proposition \ref{pra1} gives also some new necessary conditions on the invariant measure to be spectral, in the case of equal probabilities.

We denote by $\hat\mu$, the Fourier transform of the measure $\mu$:
$$\hat\mu(t)=\int e^{2\pi it\cdot x}\,d\mu(x),\quad(t\in\br^n).$$

\begin{proposition}\label{pra1}
Let $\mu_B$ be the invariant measure for the affine IFS $(\tau_b)_{b\in B}$ (equal probabilities).
Suppose that $\Lambda$ is a spectrum for $\mu_B$. Then
\begin{equation}
	\sum_{\lambda\in\Lambda}|\hat\mu_B(t-S^{-1}\lambda)|^2=N,\quad(t\in\br^n),
	\label{eqa1_1}
\end{equation}

\begin{equation}
	\sum_{\lambda\in\Lambda}e^{-2\pi iS^{-1}\lambda\cdot(b-b')}e^{2\pi i( t\cdot b-t'\cdot b')}\hat\mu_B(t-S^{-1}\lambda)\cj{\hat\mu_B(t'-S^{-1}\lambda)}=0,\quad(b\neq b').
	\label{eqa1_2}
\end{equation}
\end{proposition}

\begin{proof}
By Theorem \ref{thnoo} there is no overlap.

Since there is no overlap, we have $\mu_B(\tau_b(X_B))=\frac1N$ for all $b\in B$. We compute the Fourier coefficients of the function 
$f_{t,b}=e_t\chi_{\tau_b(X_B)}$. 

$$\ip{f_{t,b}}{e_\lambda}=\int e^{2\pi i t\cdot x}\chi_{\tau_b(X_B)}e^{-2\pi i \lambda\cdot x}\,d\mu(x)=\frac1N\sum_{b'\in B}\int e^{2\pi i(t-\lambda)\cdot R^{-1}(x+b)}\chi_{\tau_b(X_B)}(\tau_{b'}(x))\,d\mu(x)=$$
since there is no overlap only one term in the sum remains, $\chi_{\tau_b(X_B)}(\tau_{b'}(x))=\delta_{b,b'}$ for all $x$,
$$=\frac{1}{N}e^{2\pi i(t-\lambda)\cdot R^{-1}b}\int e^{2\pi iS^{-1}(t-\lambda)\cdot x}\,d\mu=\frac1Ne^{2\pi iS^{-1}(t-\lambda)\cdot b}\hat\mu(S^{-1}(t-\lambda)).$$

Applying the Parseval relation and changing the variable $S^{-1}t\mapsto t$, we obtain 
$$\frac{1}{N}=\|f_{t,b}\|^2=\sum_{\lambda\in\Lambda}\frac{1}{N^2}|\hat\mu(S^{-1}(t-\lambda))|^2,$$
and this implies \eqref{eqa1_1}.

Since there is no overlap, the functions $f_{t,b}, f_{t',b'}$ are orthogonal if $b\neq b'$. Applying the Parseval relation again we obtain
\eqref{eqa1_2}
\end{proof}

Next we show that under some extra assumptions we do get that $B$ is a spectral set, as in part (c) of the \L aba-Wang Conjecture.

{\bf Notations:} 
\begin{equation}
	\delta_B=\frac1N\sum_{b\in B}\delta_b,\quad \hat\delta_B(x)=\frac{1}{N}\sum_{b\in B}e^{2\pi ib\cdot x}.
	\label{eqdb}
\end{equation}

$$\operatorname*{Per}(\hat\delta_B):=\{p\in\br^n\,|\, \hat\delta_B(x+p)=\hat\delta_B(x),\mbox{ for all }x\in\br^n\}.$$

\begin{proposition}\label{pr4.8}
Let $\mu=\mu_B$ be the invariant measure for the IFS $(\tau_b)_{b\in B}$ (equal probabilities). Suppose $\mu$ has a spectrum $\Lambda$ with the following property: there exist $a_1,\dots,a_p\in\br^n$, and $\Lambda_1,\dots,\Lambda_p\subset\operatorname*{Per}(\hat\delta_B)$ such that 
$$\Lambda=\bigcup_{i=1}^p(a_i+S\Lambda_i),\mbox{ disjoint union},$$
and all $\Lambda_i$ are spectra for $\mu$.

Then $p=N$ and $(\delta_B,\{S^{-1}a_1,\dots,S^{-1}a_p\})$ is a spectral pair.
\end{proposition}

\begin{proof}
Taking the Fourier transform of the invariance equation we get (see e.g. \cite{DJ07a,DJ07b}):
$$\hat\mu(x)=\hat\delta_B(S^{-1}x)\hat\mu(S^{-1}x)$$Since $\Lambda$ and all $\Lambda_i$ are spectra, we have for all $t\in\br^n$:
$$1=\sum_{\lambda\in\Lambda}|\hat\mu(St-\lambda)|^2=\sum_{\lambda\in\Lambda}|\hat\delta_B(t-S^{-1}\lambda)|^2|\hat\mu(t-S^{-1}\lambda)|^2$$
$$=\sum_{i=1}^p|\hat\delta_B(t-S^{-1}a_i)|^2\sum_{\lambda\in\Lambda_i}|\hat\mu(t-S^{-1}a_i-\lambda)|^2=\sum_{i=1}^p|\hat\delta_B(t-S^{-1}a_i)|^2.$$
This implies that $\{S^{-1}a_i\,|\,i\in\{1,\dots,p\}\}$ is a spectrum for the measure $\delta_B$ (see e.g. \cite{LW02,DJ07b}). Therefore $p=N$, and $(\delta_B,\{S^{-1}a_1,\dots,S^{-1}a_p\})$ is a spectral pair.

\end{proof}

In view of Propostion \ref{pr4.8} and the counterexample mentioned in Remark \ref{rem4.7}, we reformulate part (c) of The \L aba-Wang conjecture:

\begin{conjecture}\label{co4.9}
Suppose the dimension is $n=1$ and let $\mu$ be the invariant measure of the affine IFS $(\tau_b)_{b\in B}$. Assume the measure $\mu$ is spectral. Then $B$ is a spectral set.
\end{conjecture}

The measures $\mu$ on the real line of the Bernoulli class arise as infinite convolutions, see e.g., \cite{JKS08,HL08}. These infinite convolutions are widely studied and were first considered by Erd\" os.  The measures are specified by a scaling parameter $\rho$; see Proposition \ref{pr4_9}. So we have a one-parameter family of measures $\mu,\rho$, with $\rho = \frac12$ the Lebesgue measure.

In \cite{JKS08}, the authors address the particular case $\rho =3/4$, and they ask whether $\mu_{3/4}$  is spectral or not, and they offer negative evidence. Our present result settles the question since $\mu_{3/4}$ is ``overlap''.  In \cite{JKS08} the authors  display explicit infinite and maximal orthogonal families of $e_\lambda$'s  in $L^2(\mu_{3/4})$.

\begin{proposition}\label{pr4_9}
Consider the Bernoulli convolution $\mu_\lambda$, i.e., the invariant measure for the IFS $\tau_+(x)=\lambda x+1$, $\tau_-(x)=\lambda x-1$, with $\lambda\in(\frac12,1)$. Then $\mu_\lambda$ is not a spectral measure.
\end{proposition}

\begin{proof}
Suppose $\mu_\lambda$ is spectral. 
The attractor of the IFS is $X=\left[-\frac{1}{1-\lambda},\frac{1}{1-\lambda}\right]$. We know that the support of the invariant measure of an IFS is the attractor, in our case $\supp{\mu_\lambda}=X$. A simple computation shows that $\tau_+(X)$ intersects $\tau_-(X)$ in a proper interval. Since the support of the $\mu_\lambda$ is the entire interval $X$, it follows that there is overlap. Therefore, by Theorem \ref{thnoo}, the measure cannot be spectral.
\end{proof}

We can generalize Proposition \ref{pr4_9}:

\begin{proposition}
Suppose $\mu$ is a spectral measure and the support of $\mu$ is a finite union of closed intervals. Then $\mu$ is the restriction of the Lebesgue to $\supp\mu$, renormalized.
\end{proposition}

\begin{proof}
Using Proposition \ref{prjp}, intervals of the same length contained in the support of $\mu$ have the same measure. This implies that $\mu$ is proportional to the Lebesgue measure on $\supp\mu$.
\end{proof}

Another property of spectral measures is that they cannot have atoms, unless they are purely atomic and all atoms have the same measure. The result can also be found in \cite{LW06}. We include the proof for the convenience of the reader, and we give an alternative argument.

\begin{proposition}\cite{LW06}\label{prnoa}
Let $\mu$ be a spectral measure. Then $\mu$ has an atom iff $\Lambda$ is finite. In this case $\mu$ is purely atomic and all atoms have the same measure. 
\end{proposition}
\begin{proof}
Let $\Lambda$ be a spectrum for $\mu$. Suppose $a$ is an atom. Then 
$$\mu(\{a\})=\|\chi_{\{a\}}\|^2=\sum_{\lambda\in\Lambda}|\ip{e_\lambda}{\chi_{\{a\}}}|^2=\sum_{\lambda\in\Lambda}|e_\lambda(a)\mu(\{a\})|^2=\#\Lambda\cdot\mu(\{a\})^2.$$
Therefore $\Lambda$ is finite. Thus $\mu(\{a\})=\frac{1}{\#\Lambda}$.

Conversely, if $\Lambda$ is finite then $L^2(\mu)$ is finite dimensional, and this implies that $\mu$ is purely atomic. 

Another proof can be obtained from the local translation invariance: by Proposition \ref{prjp}, if there is an atom, then by translation, all points are atoms, and since the measure is finite, there can be only finitely many of them. 
\end{proof}

\section{Embeddings into the $L^2$-space of the compact Bohr group}\label{sec4}

The theory of almost periodic functions, as initially envisioned by
Harald Bohr and Abram Besicovitch (the $L^2$-theory) had as its motivation in
classical and down-to earth questions from number theory and astronomy
\cite{FJ54}. Since then these ideas have found formulations in the context of
duality for locally compact abelian groups. Since we consider here IFSs and
associated probability measures $\mu$ in $\br^n$, it is natural for us to explore
problems regarding orthogonal exponentials (in $L^2(\mu)$) and spectral
duality within the model suggested by Bohr and Besicovitch. We do this in
detail below, allowing the Bohr-compactification $G$ of $\br^n$ as a universal
``receptor'' of spectral models in the context of $L^2(\br^n, \mu)$. Some
advantages of this is that a variety of disparate questions about spectral
duality acquire a unified framework, and we are able to rely on results from
the theory of compactifications.

\begin{definition}\label{def5.1}
For $(\br^n,+)$ consider the unitary characters $e_\lambda(x):=e^{2\pi i\lambda\cdot x}, {x\in\br^n}$. By duality, we adapt the terminology
\begin{equation}
	\ip{\lambda}{x}=e^{2\pi i\lambda\cdot x}=\ip{x}{\lambda}\mbox{ for }\lambda,x\in\br^n.
	\label{eq5.1}
\end{equation}
We shall use Pontryagin's duality for locally compact abelian groups $H$, i.e., 
\begin{equation}
	\widehat H:=\left\{\chi:H\rightarrow\bt\,|\,\chi\mbox{ continuous and } \chi(h_1+h_2)=\chi(h_1)\chi(h_2), \chi(-h)=\cj{\chi(h)}, \, h,h_1,h_2\in H\right\}.
	\label{eq5.2}
\end{equation}
($\bt:=\{z\in\bc\,|\,|z|=1\}$)
Moreover, $\widehat H$ is given the compact-open topology. 

It is known that $\widehat H$ is again a locally compact abelian group under the operation $(\chi_1\chi_2)(h):=\chi_1(h)\chi_2(h)$, $h\in H$. 
Moreover:

(i) $\widehat{\widehat H}\cong H$, i.e., the natural embedding $H\hookrightarrow\widehat{\widehat H}$ is onto.

(ii) $H$ is compact iff $\widehat H$ is discrete.

We apply  this to $(\br^n,+)$ where $\br^n$ is given the usual topology. When it is equipped with the discrete topology, it is denoted $\br^n_{disc}$

(iii) It follows that $G:=\widehat{\br^n_{disc}}$ is a compact abelian group, by (i)-(ii), with normalized Haar measure $\mu_{Bohr}$, where the subscript "`Bohr"' is after Harald Bohr, \cite{BeBo31,BS32,FJ54}.

(iv) Dualizing the natural mapping $\br^n_{disc}\hookrightarrow\br^n$ (continuous!) we get 
\begin{equation}\label{eq5.4}
\br^n=\widehat{\br^n}\hookrightarrow G;
\end{equation}
 i.e., $\br^n$ is naturally embedded into $G$: hence the name ``{\it Bohr compactification}''.

For $T>0$ set $Q_T=\{x\in\br^n\,|\, -T\leq x_j\leq T, 1\leq j\leq n\}$. Bohr proved that the following limit exists for the almost periodic functions
\begin{equation}
	\lim_{T\rightarrow\infty}\frac{1}{(2T)^n}\int_{Q_T}f(x)\,dx=:\lim_{T\rightarrow\infty}\langle f\rangle_T=\langle f\rangle.
	\label{eq5.5}
\end{equation}
\end{definition}

\begin{definition}\label{def5.2}
A continuous function $f$ on $\br^n$ is said to be {\it almost periodic} if for all $\epsilon>0$ there exists $T\in\br_+$ such that for all $y\in \br^n$ there exists $p\in y+Q_T$ such that 
\begin{equation}
	|f(x)-f(x+p)|<\epsilon, \mbox{ for all } x\in\br^n.
	\label{eq5.6}
\end{equation}
Moreover, if $f$ is almost periodic, then 
\begin{equation}
	\lim_{T\rightarrow\infty}\langle f\rangle_T=\langle f\rangle=\int_G f\,d\mu_{Bohr},
	\label{eq5.7}
\end{equation}
where $\br^n$ is embedded in $G$ via \eqref{eq5.4}. In particular, a continuous almost periodic function $\br^n$ extends naturally to a continuous function on $G$. 
\end{definition}

We now get the following

\begin{theorem}\label{th5.3} Abstract $L^2$-embedding. Let $\xi:\br^n\rightarrow G$ denote Bohr's embedding \eqref{eq5.4}, i.e., 
\begin{equation}
	\ip{\xi(x)}{\lambda}=e_\lambda(x)=:e^{2\pi i\lambda\cdot x}, \quad (x,\lambda\in\br^n)
	\label{eq5.8}
\end{equation}
and set
\begin{equation}
	\tilde e_\lambda(\chi):=\chi(\lambda),\quad(\lambda\in\br^n_{disc},\chi\in G)
	\label{eq5.9}
\end{equation}

Let $\mu$ be a finite measure on $\br^n$ and let $\Lambda\subset\br^n$ be the subset of $\br^n$. Then the set $E(\Lambda):=\{e_\lambda\}_{\lambda\in\Lambda}$ is orthonormal in $L^2(\mu)$ iff the embedding given by \eqref{eq5.8}-\eqref{eq5.9}
\begin{equation}
	W_\Lambda:e_\lambda\mapsto \tilde e_\lambda\mbox{ where }\tilde e_\lambda(\chi)=\chi(\lambda),\, \chi\in G, W_\Lambda:\H_\Lambda:=\operatorname*{ cl span}\{e_\lambda\}\hookrightarrow L^2(G).
	\label{eq5.10}
\end{equation}
is an isometric operator.
\end{theorem}

\begin{remark}
In the hypotheses of Theorem \ref{th5.3}, we have the following inclusions: $$\Lambda\subset\br^n_{disc}\hookrightarrow\br^n$$ The dual of this inclusion is $$G=\widehat{\br^n_{disc}}\hookleftarrow\widehat\br^n=\br^n.$$ The isometry $W_\Lambda$ maps $$E(\Lambda)\subset L^2(\br^n,\mu)\rightarrow L^2(G)$$

\end{remark}
\begin{proof}[Proof of Theorem \ref{th5.3}]
We only need to check that if $\lambda$ and $\lambda'$ are distinct points in $\Lambda$, then 
\begin{equation}
	\ip{e_\lambda}{e_{\lambda'}}_{L^2(G)}=0.
	\label{eq5.11}
\end{equation}
But we may compute \eqref{eq5.11} with the use of \eqref{eq5.7}:

\begin{equation}
	\ip{e_\lambda}{e_{\lambda'}}_{L^2(G)}=\lim_{T\rightarrow\infty}\frac{1}{(2T)^n}\int_{Q_T} e^{2\pi i(\lambda-\lambda')\cdot x}\,dx=0
	\label{eq5.12}
\end{equation}
by a direct computation.

\end{proof}

%
%
%
%
%

\begin{lemma}\label{lem5.8}
The family $\{e_\lambda\}_{\lambda\in\br^n}$ is linearly independent in $L^2(\mu)$ if $\mu$ is spectral and has an infinite spectrum.  
\end{lemma}

\begin{proof}
If $\mu$ has an infinite spectrum then the measure $\mu$ has no atoms, by Proposition \ref{prnoa}. 
Suppose 
\begin{equation}
	\sum_{\mbox{ finite,} \lambda\in\br^n}\xi_\lambda e_\lambda=0\mbox{ in } L^2(\mu).
	\label{eq5.16}
\end{equation}
The sum in \eqref{eq5.16} is a trigonometric polynomial so it can only have a discrete at most countable set of zeros. From \eqref{eq5.16} we conclude that the measure $\mu$ is supported on this set. Therefore it has to be an atomic measure. The contradiction implies the lemma.
\end{proof}

\begin{theorem}\label{th5.12}
Let $\mu$ be a probability measure on $\br^n$ and $\Lambda\subset\br^n$. The following are equivalent
\begin{enumerate}
\item
The set $E(\Lambda):=\{e_\lambda\,|\,\lambda\in\Lambda\}$ is orthogonal in $L^2(\mu)$.
\item
Let $\H(\Lambda):=\operatorname*{clspan}_{L^2(\mu)|} E(\Lambda)$. The operator
$W_\Lambda:\H(\Lambda)\rightarrow L^2(G)$, $W_\Lambda(e_\lambda)=\tilde e_\lambda$ for $\lambda\in\Lambda$, is isometric. 
\item 
The function
\begin{equation}
h_\Lambda(t):=\sum_{\lambda\in\Lambda}|\hat\mu(t-\lambda)|^2	
	\label{eq5.12.1}
\end{equation}

satisfies the inequality $h_\Lambda(t)\leq 1$ for all $t\in\br^n$.

\end{enumerate}
\end{theorem}

\begin{proof}
We already compared (i)$\Leftrightarrow$(ii). 

(i)$\Rightarrow$(iii). Using Bessel's inequality we have
$$1=\|e_t\|_{L^2(\mu)}^2\geq \sum_{\lambda\in\Lambda}|\ip{e_\lambda}{e_t}|^2=\sum_{\lambda\in\Lambda}|\hat\mu(t-\lambda)|^2=h_\Lambda(t),$$
for all $t\in\br^n$. So (iii) holds.

(iii)$\Rightarrow$(i). For any $\lambda_0\in\Lambda$ we have 
$$1\geq h_{\Lambda}(\lambda_0)=|\hat\mu(\lambda_0-\lambda_0)|^2+\sum_{\lambda\in\Lambda,\lambda\neq\lambda_0}|\hat\mu(\lambda_0-\lambda)|^2=1+\sum_{\lambda\in\Lambda,\lambda\neq\lambda_0}|\hat\mu(\lambda_0-\lambda)|^2.$$
Therefore $\hat\mu(\lambda_0-\lambda)=0$ for all $\lambda\neq \lambda_0$. But this implies that $\ip{e_\lambda}{e_{\lambda_0}}=0$, so (i) holds.
\end{proof}

\begin{remark}
It is known that a set $\Lambda$ forms a spectrum for a measure iff $h_\Lambda$ is constant 1 (see e.g. \cite{JP98b,DJ07b}). We note here that by Theorem \ref{th5.12} the orthogonality (without completeness) can be characterized in terms of $h_\Lambda$. Moreover the maximal orthogoanlity, within the class of exponential functions can be expressed in terms of $h_\Lambda$, as we prove below in Theorem \ref{th5.13}.
\end{remark}

\begin{theorem}\label{th5.13}
Let $\mu$ be a probability measure on $\br^n$ and $\Lambda\subset\br^n$. Then $E(\Lambda)$ is a maximal family of orthogonal exponentials if and only if  
$0<h_\Lambda(t)\leq 1$ for all $t\in\br^n$.
\end{theorem}

\begin{proof}
If $E(\Lambda)$ is maximal orthogonal, then we know $h_\Lambda(t)\leq 1$ from Theorem \ref{th5.12}. Suppose there is some $t_0$ with $h_\Lambda(t_0)=0$. Then 
$\hat\mu(t_0-\lambda)=0$ for all $\lambda\in\Lambda$ so $\ip{e_\lambda}{e_{t_0}}=0$ for all $\lambda$. Since $E(\Lambda)$ is maximal this implies that $t_0\in\Lambda$, but then $1=\mu(t_0-t_0)=0$, a contradiction.

Conversely, if $0<h_\Lambda \leq 1$, then,  from Theorem \ref{th5.12} we know that $E(\Lambda)$ is orthogonal. Suppose it is not maximal, so there is $t_0\not\in\Lambda$ such that $e_{t_0}\perp e_\lambda$ for all $\lambda\in\Lambda$. Then $\hat\mu(t_0-\lambda)=0$ so $h_\Lambda(t_0)=0$, a contradiction.
\end{proof}

Next, we show that the embedding of a spectral measure into the Bohr group intertwines the local translations. We recall first the definition of translations on the Bohr group.
\begin{lemma}\label{lem4.1.1}
Let $G=\widehat{\br^n_{disc}}$ be the Bohr group. For $a \in \br^n$ and $\chi\in G$, set
\begin{equation}
	(a\cdot\chi)(x):=e^{2\pi i a\cdot x}\chi(x),\quad(x\in\br^n).
	\label{eq4.1.1}
\end{equation}
The mapping 
\begin{equation}
\br^n\times G\ni(a,\chi)\mapsto a\cdot \chi\in G	
	\label{eq4.1.2}
\end{equation}
is a continuous transformation group, i.e., 
\begin{equation}
	(a+b)\cdot\chi=a\cdot(b\cdot\chi),
	\label{eq4.1.3}
\end{equation}
holds for all $a,b\in\br^n$ and $\chi\in G$.
\end{lemma}

\begin{proof}
The continuity assertion is clear from the definition of the topology on $G$ (i.e., generated by the cylinder set neighborhoods).

For the verification of \eqref{eq4.1.3}, let $a,b,x\in\br^n$, and $\chi\in G$. Then
$$((a+b)\cdot\chi)(x)=e_{a+b}(x)\chi(x)=e_a(x)e_b(x)\chi(x)=e_a(x)(b\cdot\chi)(x)=(a\cdot(b\cdot\chi))(x)$$
i.e., the desired formula \eqref{eq4.1.3} holds.
\end{proof}

\begin{corollary}\label{cor4.1.2}
Let $G=\widehat{\br^n_{disc}}$ be the Bohr group with Haar measure $\mu_{Bohr}$. Then there is a natural strongly continuous unitary representation $\U_{Bohr}$ of $\br^n$ acting on $L^2(G,\mu_{Bohr})$ by
\begin{equation}
	(\U_{Bohr}(a)f)(\chi):=f(a\cdot\chi),\quad(a\in\br^n,\chi\in G, f\in L^2(G,\mu_{Bohr}))
	\label{eq4.1.4}
\end{equation}
\end{corollary}

\begin{proof}
The fact that \eqref{eq4.1.4} defines a strongly continuous representation of $(\br^n,+)$ acting on $L^2(G,\mu_{Bohr})$ is immediate from Lemma \ref{lem4.1.1}
\end{proof}

\begin{theorem}\label{th4.1.4}
Let $\mu$ be a Borel probability measure on $\br^n$. Suppose $\mu$ is spectral with spectrum $\Lambda$. Let $\xi:\br^n\hookrightarrow G$ be the embedding into the Bohr group, see Definitions \ref{def5.1}, \ref{def5.2}. Let $U_\Lambda=:U$ be the unitary representation of $\br^n$ on $L^2(\mu)$ by local translations as in Definition \ref{defft}.

Let 
\begin{equation}
	W:L^2(\mu)\rightarrow L^2(G)
	\label{eq4.1.6}
\end{equation}
be the isomorphic embedding from Theorem \ref{th5.3}. Then the following intertwining relation holds:
\begin{equation}
	\U_{Bohr}(a)W=WU_\Lambda(a),\quad(a\in\br^n)
	\label{eq4.1.7}
\end{equation}
\end{theorem}

\begin{proof}
Since $\{e_\lambda\}_{\lambda\in\Lambda}$ is an ONB in $L^2(\mu)$, we only need to verify the operator commutation relation \eqref{eq4.1.7} on these basis vectors.

Indeed, let $a\in\br^n$, $\lambda\in\Lambda$ and $\chi\in G$ be given. Then we get \eqref{eq4.1.7} by the following computation:
$$((\U_{Bohr}(a)W)e_\lambda)(\chi)=(We_\lambda)(a\cdot\chi)\stackrel{\mbox{by \eqref{eq5.8}-\eqref{eq5.10}}}{=}(a\cdot\chi)(\lambda)=e_\lambda(a)\chi(\lambda)=e_\lambda(a)\tilde e_\lambda(\chi)=(W\U_\Lambda(a)e_\lambda)(\chi).$$

Since this holds for all $\chi\in G$, the desired formula \eqref{eq4.1.7} follows. 
\end{proof}

\begin{corollary}\label{cor4.1.4}

Let $\mu,\Lambda$ and $G$ be as specified in Theorem \ref{th4.1.4}, let $x,a\in\br^n$ be chosen such that $x\in\supp\mu$, and $x+a\in\supp\mu$. Let $f\in L^2(\mu)$. Then 
\begin{equation}
	(U_\Lambda(a)f)(x)=f(x+a)
	\label{eq4.1.8}
\end{equation}
\end{corollary}
\begin{proof}
The result follows from Theorem \ref{th4.1.4} since $f$ can be expanded in the ONB $(e_\lambda)_{\lambda\in\Lambda}$ in $L^2(\mu)$; and Theorem \ref{th4.1.4} states that \eqref{eq4.1.8} holds if $f=e_\lambda$.
\end{proof}
\begin{remark}\label{rem5.8}
Consider the Bernoulli convolutions $L^2(\mu_\lambda)$ for $\lambda=1/4$ and $\lambda=3/4$; see Proposition \ref{pr4_9}. An inspection of the formula (see e.g. \cite{JKS08})
$$\hat\mu_\lambda(t)=\prod_{k=1}^\infty\cos(2\pi\lambda^k t)$$
shows that $Z(\hat\mu_{1/4})\subset Z(\hat\mu_{3/4})$ where $Z$ denotes the ``zero-set''. 

Since 
$$\Gamma:=\left\{\sum_{finite} a_i4^i\,|\,a_i\in\{0,1\}\right\}=\{0,1,4,5,16,17,\dots\}$$
makes $(\mu_{1/4},\Gamma)$ into a spectral pair, it follows that $\{e_\gamma\,|\,\gamma\in\Gamma\}\subset L^2(\mu_{3/4})$ is an infinite orthogonal set of Fourier frequencies. As a result the isometry $W: e_\gamma\in L^2(\mu_{1/4})\mapsto e_\gamma\in L^2(\mu_{3/4})$ extends by linearity to an isometric embedding of the Hilbert space $L^2(\mu_{1/4})$ into $L^2(\mu_{3/4})$. Hence the two Hilbert spaces $L^2(\mu_{1/4})$ and $WL^2(\mu_{1/4})\subset L^2(\mu_{3/4})$ have the same representation in $L^2(G)$; see Theorem \ref{th5.3}.
\end{remark}

\begin{acknowledgements}
The authors are pleased to thank the following for helpful conversations: Ka-Sing Lau, Karen Shuman, Keri Kornelson, Myung-Sin Song. We thank the anonymous referee for suggestions that improved the paper.
\end{acknowledgements}

\bibliographystyle{alpha}	
\bibliography{ffin}

\end{document}